\documentclass[reqno]{amsart}
\usepackage{amssymb,lattice}

\theoremstyle{plain}
\newtheorem{lemma}{Lemma}[section]
\newtheorem{theorem}[lemma]{Theorem}

\newtheorem{proposition}[lemma]{Proposition}
\newtheorem{corollary}[lemma]{Corollary}
\newtheorem{example}[lemma]{Example}
\newtheorem{claim}{Claim}

\newtheorem*{stat}{\name}
\newcommand{\name}{testing}

\theoremstyle{definition}
\newtheorem{definition}[lemma]{Definition}
\newtheorem{problem}{Problem}
\newtheorem{case}{Case}

\theoremstyle{remark}

\newenvironment{all}[1]{\renewcommand{\name}{#1}\begin{stat}}
{\end{stat}}

\newcommand{\pp}{\partial}
\newcommand{\ppu}{\partial_{\mathrm{u}}}
\newcommand{\ppm}{\partial_{\mathrm{m}}}

\newcommand{\tr}{\vartriangleleft}

\newcommand{\utr}{\trianglelefteq}

\newcommand{\dtr}{\mathbin{\vartriangleleft\kern-10pt {\lower
3pt\hbox{$\scriptscriptstyle\ne$}}\kern3pt}}

\newcommand{\DD}{\mathbin{\gd}}

\DeclareMathOperator{\Con}{Con}

\DeclareMathOperator{\J}{J}
\DeclareMathOperator{\JE}{JE}
\DeclareMathOperator{\je}{je}
\DeclareMathOperator{\HH}{H}
\DeclareMathOperator{\At}{At}

\author{G.~Gr\"atzer}
\thanks{The research of the first author was  supported by the
NSERC of Canada.}
\address{Department of Mathematics\\ University of Manitoba\\
Winnipeg MN, R3T 2N2\\ Canada}
\email{gratzer@cc.umanitoba.ca}
\urladdr{http://server.maths.umanitoba.ca/homepages/gratzer/}

\author{F.~Wehrung}
\address{CNRS ESA 6081\\
Universit\'e de Caen, Campus II\\
D\'epartement de Math\'ematiques\\
BP 5186\\
14032 Caen Cedex\\
France}
\email{wehrung@math.unicaen.fr}
\urladdr{http://www.math.unicaen.fr/\~{}wehrung}


\begin{document}

\title[Number of join-irreducibles in a congruence
representation]%
{On the number of join-irreducibles in a
congruence representation of a finite distributive lattice}

\begin{abstract} For a finite lattice $L$, let $\utr_L$ denote
the reflexive and transitive closure of the join-dependency
relation on $L$, defined on the set $\J(L)$ of all
join-irreducible elements of $L$. We characterize the relations
of the form $\utr_L$, as follows:

\begin{all}{Theorem} Let $\utr$ be a quasi-ordering on a finite
set $P$. Then the following conditions are equivalent:
\begin{enumerate}
\item[\tup{(i)}] There exists a finite lattice $L$ such that
$\vv<\J(L),\utr_L>$ is isomorphic to the quasi-ordered set
$\vv<P,\utr>$.
\item[\tup{(ii)}] $|\setm{x\in P}{p\utr x}| \neq 2$, for any
$p\in P$.
\end{enumerate}
\end{all}

For a finite lattice $L$, let $\je(L) = |\J(L)| - |\J(\Con L)|$,
where $\Con L$ is the congruence lattice of $L$. It is
well-known that the inequality $\je(L)
\geq 0$ holds. For a finite distributive lattice $D$, let us
define the
\emph{join-excess function}:
\[
\JE(D) = \min(\je(L) \mid \Con L \iso D). \] We provide a
formula for computing the join-excess function of a finite
distributive lattice $D$. This formula implies that
$\JE(D)\leq(2/3)|\J(D)|$, for any finite distributive lattice
$D$; the constant $2/3$ is best possible.

A special case of this formula gives a characterization of
congruence lattices of finite lower bounded lattices.
\end{abstract}

\maketitle

\section*{Introduction}

In \cite{GS62}, the first author and E.\,T. Schmidt proved the
following result:

\begin{all}{Representation by finite sectionally complemented
lattices} For every finite distributive lattice $D$, there
exists a finite, sectionally complemented lattice~$L$ such that
the congruence lattice $\Con L$ of $L$ is isomorphic to $D$.
Furthermore, if $n$ denotes the number of join-irreducible
elements of $D$ and if $n>0$, then $L$ has fewer than $2n$
join-irreducible elements.
\end{all}

On the other hand, it follows from classical results of lattice
theory, see Theorem~\ref{T:DCon}, that for any finite lattice
$L$, if $\J(L)$ denotes the set of all join-irreducible elements
of $L$, then the inequality
\[ |\J(L)|\geq|\J(\Con L)|
\] holds. So if we define $\je(L) = |\J(L)| - |\J(\Con L)|$,
then $\je(L) \geq 0$ and $\je(L)$ is one measure of the
efficiency of the representation of $D =
\Con L$ as a congruence lattice. Define
\[
\JE(D) = \min(\je(L) \mid \Con L \iso D). \] Then from this
point of view, the best representation of a finite distributive
lattice~$D$ as a congruence lattice of a finite lattice $L$ is
obtained when $\je(L) = \JE(\Con L)$.

If $D$ is a finite distributive lattice with $n$
join-irreducible elements and if $n>0$, then the least number
$\JE(D)$ satisfies the inequality
\begin{equation}\label{Eq:n2n} 0 \leq \JE(D) < n.
\end{equation} In this paper, we shall give a formula that
computes $\JE(D)$ from $D$, see Theorems~\ref{T:Ineq} and
\ref{T:OptJoin}. We would like to emphasize that our formula
does not estimate $\JE(D)$ but gives the
\emph{exact} value. However, it implies the better estimate
\[ 0 \leq \JE(D) \leq \frac{2}{3}n,
\] and the constant $2/3$ in this estimate is best possible (see
Corollary~\ref{C:5/3bound}).

The formula that computes $\JE(D)$ from $D$ is extremely
``fast'' (in linear time) and it only involves properties of the
``upper layer'' of $\J(D)$---more precisely, the maximal
elements of $\J(D)$ and the elements that they cover.

The basic new concept is a \emph{spike}. A spike of a finite
poset $P$ is a pair
$\vv<p,q>$ of elements of $P$ such that $q$ is maximal and $q$
is the unique element of $P$ that covers~$p$.

As a corollary, we obtain a characterization of those $D$ that
are isomorphic to $\Con L$ for some finite lattice $L$ that is
\emph{lower bounded} in the sense of R.\,N. McKenzie~\cite{rM70}
(see also  R.~Freese, J.~{}Je\v zek, and J.\,B. Nation
\cite{FJN95}), or equivalently, \emph{amenable} in the sense of
our papers \cite{GW0,GW1}, see Section~\ref{S:MinNo}.

This characterization is more conveniently expressed in terms of
the poset $P=\J(D)$ of all join-irreducible elements of~$D$:
\begin{equation}\label{Eq:spfree}
\text{$P$ has no spikes}.
\end{equation} See Corollary~\ref{C:SpFree}.

We obtain these results by studying the \emph{join-dependency
relation},
$\DD_L$, on a finite lattice $L$, or, rather, its reflexive,
transitive closure, that we denote by $\utr_L$. It is well-known
that $\utr_L$ determines the congruence structure of $L$, see
Theorem~\ref{T:DCon}. Our main result, Theorem~\ref{T:Main},
describes when a binary relation on a finite set is isomorphic
to $\utr_L$ on $\J(L)$, for some finite lattice $L$. This
description is very similar to condition \eqref{Eq:spfree}.

Another consequence of Theorem~\ref{T:Main} is the
characterization, for a finite lattice~$L$, of the canonical
surjective map from $\J(L)$ onto $\J(\Con L)$, see
Theorem~\ref{T:Part}.

In all these results the finite lattice $L$ we construct is
\emph{atomistic}, that is, every element is a (finite) supremum
of atoms. This is not surprising, in view of the result of M.
Tischendorf \cite{Tisch}: Every finite lattice $K$ has a finite,
atomistic, congruence-preserving extension $L$; in addition,
$\vv<\J(K),\utr_K>\iso\vv<\J(L),\utr_L>$.

By G. Gr\"atzer and E.\,T. Schmidt \cite{GS}, every finite
lattice has a finite, \emph{sectionally complemented},
congruence-preserving extension. However, in our results $L$
cannot be taken as sectionally complemented. In
Example~\ref{Ex:nosc}, we describe a finite distributive lattice
$D$ that can be represented as $\Con L$ for $L$ finite, lower
bounded, atomistic, but which cannot be represented as $\Con L$
for $L$ finite, lower bounded, sectionally complemented.

\section{Basic concepts}

Let $L$ be a finite lattice. We denote by $\J(L)$ the set of all
join-irreducible elements of $L$. For $p\in L$, we denote by
$p_*$ the unique element of $L$ covered by $p$. The
\emph{join-dependency relation}, $\DD_L$, is the binary relation
defined on $\J(L)$ by
\[ p\DD_Lq\quad\text{if{}f}\quad p\ne q,\ \text{and there
exists }  x\in L\text{ such that } p\leq q\jj x\text{ and
}p\nleq q_*\jj x. \] In particular, note that $p\DD_Lq$ implies
that $p\nleq q$.

A useful alternative description of the join-dependency relation
on $\J(L)$ arises from \emph{minimal pairs}. Let $L$ be a finite
lattice. For $J$, $I \ci L$, we say that $I$
\emph{dominates} $J$, in notation, $J\ll I$, if{}f for all
$x\in J$, there exists $y\in I$ such that $x\leq y$. As in
H.\,S. Gaskill
\cite{hG73}, and H.\,S. Gaskill, G. Gr\"atzer, and C.\,R. Platt
\cite{GGP75}, a minimal pair of a finite lattice $L$ is a pair
$\vv<p,I>$, where
$p\in\J(L)$, $I\ci\J(L)$, $p\nin I$, $p\leq\JJ I$, and, for
every subset $J$ of $\J(L)$ such that $J\ll I$, the inequality
$p\leq\JJ J$ implies that $I\ci J$. Observe that if $\vv<p,I>$
is a minimal pair, then $I$ has at least two elements.

\begin{lemma}\label{L:DescrD} Let $L$ be a finite lattice. For
all $p$,
$q\in\J(L)$, the following are equivalent:
\begin{enumerate}
\item[\tup{(i)}] $p\DD_Lq$.
\item[\tup{(ii)}] There exists $I\ci\J(L)$ such that $\vv<p,I>$
is a minimal pair of $L$ and $q\in I$.
\end{enumerate}
\end{lemma}

See, for example, Lemma~2.31 in R.~Freese, J.~{}Je\v zek, and
J.\,B. Nation \cite{FJN95}.

We shall denote by $\tr_L$ (resp., $\utr_L$) the transitive
closure (resp., reflexive transitive closure) of the relation
$\DD_L$. By definition, $\utr_L$ is a quasi-ordering on $\J(L)$,
that is, it is reflexive and transitive. Moreover, we will
denote by $\asymp_L$ the equivalence relation associated with
$\utr_L$; so, for $p$, $q\in\J(L)$,
\[ p\asymp_Lq\quad\text{if{f}}\quad p\utr_Lq\utr_L p. \]
We refer to Theorem~2.30 and Lemma~2.36 in R.~Freese, J.~{}Je\v
zek, and J.\,B. Nation~\cite{FJN95} for a proof of the following
result:

\begin{theorem}\label{T:DCon} Let $L$ be a finite lattice. For
all
$p\in\J(L)$, let $\gQ(p)$ be the congruence of $L$ generated by
the pair
$\vv<p_*,p>$. Then the following statements hold:
\begin{enumerate}
\item[\tup{(i)}] $\gQ$ is a map from $\J(L)$ onto $\J(\Con L)$.
\item[\tup{(ii)}] For all $p$, $q\in\J(L)$, $\gQ(p)\ci\gQ(q)$
if{f} $p\utr_Lq$.
\end{enumerate}
\end{theorem}

We shall use the following notation. If $\tr$ is a binary
relation on a set
$P$, then we define the \emph{upper $\tr$-segment} of $p$ as
\[ [p]^{\tr}=\setm{x\in P}{p\tr x},
\] for any $p\in P$.

\section{The relations $\DD$, $\tr$, $\utr$ on a finite lattice}

The elementary properties of $\utr_L$ will be described in
Proposition~\ref{P:trutr}; to prepare for it, we first prove a
simple lemma:

\begin{lemma}\label{L:Drefl} Let $L$ be a finite lattice. For
all $p\in\J(L)$, if $p\tr_L p$, then there exists $x\in\J(L)$
such that $x\ne p$ and
$p\tr_Lx\tr_Lp$. \end{lemma}

\begin{proof} By the definition of $\DD_L$, one cannot have
$p\DD_Lp$. Therefore, by the definition of $\tr_L$, there is a
positive integer $n$ and there are elements $x_0$, \dots,
$x_n\in\J(L)$ such that
\[ p=x_0\DD_Lx_1\DD_L\cdots\DD_Lx_n=p.
\] Then $x_1\ne p$, and $p\tr_Lx_1\tr_Lp$, so that $x=x_1$
satisfies the required conditions.
\end{proof}

As a consequence, each of the relations $\tr_L$ and $\utr_L$ can
be defined in terms of the other:

\begin{proposition}\label{P:trutr} Let $L$ be a finite lattice.
For all $p$,
$q\in\J(L)$, the following statements hold:
\begin{enumerate}
\item[\tup{(i)}] $p\utr_Lq$ if{f} $p\tr_Lq$ or
$p=q$.\label{I:first}
\item[\tup{(ii)}] $p\tr_Lq$ if{f} one of the two following
conditions hold:
\begin{enumerate}
\item[\tup{(a)}] $p\utr_Lq$ and $p\ne q$. \item[\tup{(b)}]
$|\,[p]^{\asymp_L}\,|\geq2$ and $p=q$. \end{enumerate}
\end{enumerate}
\end{proposition}

\begin{proof}\hfill

(i) is trivial.

Now we prove (ii). Let us assume first that $p\tr_Lq$. If $p\ne
q$, then (a) holds. Now assume that $p=q$. By
Lemma~\ref{L:Drefl}, there exists $x\ne p$ in
$\J(L)$ such that $p\tr_Lx\tr_Lp$. In particular, $x$ belongs to
$[p]^{\asymp_L}$, so that $|\,[p]^{\asymp_L}\,| \geq 2$. Thus
(b) holds.

Conversely, (a) trivially implies that $p\tr_Lq$. Assume (b).
Since
$|[p]^{\asymp_L}|\geq2$, there exists $x\ne p$ such that
$p\asymp_Lx$. Necessarily,
$p\tr_Lx\tr_Lp$, so that $p\tr_Lp=q$.
\end{proof}

\begin{proposition}\label{P:Constr} Let $L$ be a finite lattice.
Then
$|\,[p]^{\utr_L}\,|\ne2$, for any $p\in\J(L)$.
\end{proposition}

\begin{proof} Assume that $[p]^{\utr_L}$ has exactly two
elements. In particular, there exists
$q_0\ne p$ such that $p\utr_Lq_0$. Therefore, $p\tr_Lq_0$, so
there exists $q$ such that $p\DD_Lq$ and $q\utr_Lq_0$. Since
$p\DD_Lq$, there exists, by Lemma~\ref{L:DescrD}, a subset $I$
of $\J(L)$ such that $\vv<p,I>$ is a minimal pair of $L$ and
$q\in I$. In particular, $|I|\geq2$, so $I$ contains some $x\ne
q$. Since $x\in I$, $x$ is also distinct from $p$. Therefore,
$\set{p,q,x}\ci[p]^{\utr_L}$, a contradiction.
\end{proof}

\section{Finite atomistic lattices from the $\utr$ relation} For
a finite atomistic lattice $L$, let $\At(L)$ denote the set of
atoms of $L$. Of course,
$\At(L) = \J(L)$.

If $\utr$ is a quasi-ordering on a set $P$, we denote by $\dtr$
the binary relation on $P$ defined by
\[ p\dtr q\quad\text{if{}f}\quad p\utr q\text{ and }p\ne q. \]
The main goal of this section is to prove the following converse
of Proposition~\ref{P:Constr}:

\begin{theorem}\label{T:Main} Let $P$ be a finite set, let
$\utr$ be a quasi-ordering on $P$. Then the following conditions
are equivalent:
\begin{enumerate}
\item[\tup{(i)}] There exists a finite atomistic lattice $L$
such that \[
\vv<P,\utr,\dtr>\iso\vv<\At(L),\utr_L,\DD_L>. \]
\item[\tup{(ii)}] There exists a finite lattice $L$ such that \[
\vv<P,\utr>\iso\vv<\J(L),\utr_L>.
\]
\item[\tup{(iii)}] $|\,[p]^{\utr}\,|\ne2$, for all $p\in P$.
\end{enumerate}
\end{theorem}

\begin{proof}\hfill

(i)$\Rightarrow$(ii) is trivial.

(ii)$\Rightarrow$(iii) follows from Proposition~\ref{P:Constr}.

We prove, finally, the direction (iii)$\Rightarrow$(i). So we
are given
$\vv<P,\utr>$ satisfying that $|\,[p]^{\utr}\,|\ne2$, for any
$p\in P$. Let us say that a subset $X$ of $P$ is \emph{closed},
if for all $x$,~$y\in X$ such that $x\ne y$ and
\begin{equation}\label{Eq:DefL} p\utr x,\,y\quad\text{implies
that}\quad p\in X, \end{equation} for all $p\in P$, where $p\utr
x$, $y$ stands for $p\utr x$
\emph{and} $p\utr y$. Furthermore, we denote by $L$ the set of
all closed subsets of $P$.

It is obvious that any intersection of closed subsets of $P$ is
closed, and that both $\es$ and $P$ are closed. Thus, $L$ is a
closure system in the powerset lattice of~$P$. In particular,
$L$, partially ordered by containment, is a lattice.
Furthermore, by the definition of a closed subset of $P$, it is
obvious that the singleton $\ge(p)=\set{p}$ is closed, for all
$p\in P$. Hence, the lattice $L$ is atomistic, and the atoms of
$L$ are exactly the singletons of elements of $P$. In
particular, $\ge$ is a bijection from $P$ onto $\J(L)$. We shall
now prove that $\ge$ is an isomorphism from
$\vv<P,\utr,\dtr>$ onto $\vv<\J(L),\utr_L,\DD_L>$.

For all $X\ci P$, we shall denote by $\ol{X}$ the \emph{closure}
of $X$ in
$L$, that is, the least element of $L$ that contains $X$.
\emph{A priori}, the closure of $X$ is computed by iteration of
the rule \eqref{Eq:DefL}. Our next claim will show that only one
step is required:

\setcounter{claim}{0}
\begin{claim} For every subset $X$ of $P$, the closure of $X$
can be computed by the following formula:
\begin{equation}\label{Eq:Clos}
\ol{X}=X\uu
\setm{p\in P}{p\utr x,\,y,\ \text{for some }x,\,y\text{ in }X
\text{ such that }x\ne y}.
\end{equation}
\end{claim}

\begin{proof} Let $X'$ denote the right side of \eqref{Eq:Clos}.
It is obvious that $X'$ contains $X$ and that every closed
subset of $P$ containing $X$ contains $X'$. Thus it suffices to
prove that $X'$ is closed. So let $p\in P$ and $x$, $y\in X'$
such that $x\ne y$ and $p\utr x$, $y$; we prove that $p\in X'$.
If both $x$ and $y$ already belong to $X$, then this is obvious
by the definition of $X'$. Otherwise, without loss of
generality, we can assume that
$x\in X'-X$; thus, by the definition of $X'$, there are $x_0$,
$x_1\in X$ such that $x_0\ne x_1$ and $x\utr x_0$, $x_1$. Since
$p\utr x$, it follows that
$p\utr x_0$, $x_1$; whence $p\in X'$.
\end{proof}

We conclude the proof of Theorem~\ref{T:Main} with three more
claims:

\begin{claim} The minimal pairs of $L$ are exactly the pairs of
the form
\[
\vv<\ge(p),\set{\ge(x),\ge(y)}>,
\] where $p$, $x$, $y\in P$, $p\dtr x$, $y$, and $x\ne y$.
\end{claim}

\begin{proof} Let $p$, $x$, $y$ be given as above. The join
$\ge(x)\jj\ge(y)$ of $\ge(x)$ and $\ge(y)$ in $L$ is closed and
contains $\set{x,y}$, and so it contains $p$, by the definition
of a closed subset of $P$. Hence, $\ge(p) \leq
\ge(x) \jj \ge(y)$. Since $p\nin\set{x,y}$, it follows that
$\ge(p)$ is contained neither in $\ge(x)$ nor in $\ge(y)$. Now
$\ge(x)$ and $\ge(y)$ are atoms of $L$, so
$\vv<\ge(p),\set{\ge(x),\ge(y)}>$ is a minimal pair of $L$.

Conversely, a minimal pair of $L$ has the form
$\vv<\ge(p),\ge[I]>$, where
$p\in P$, $I\ci P$, $p\nin I$, and $p\in\JJ\ge[I]$. The last
condition means that $p$ belongs to the closure of $I$, thus, by
Claim~1 and by $p\nin I$, there are $x$,
$y\in I$ such that $x\ne y$ and $p\utr x$, $y$. Let
$J=\set{x,y}$. Then
$\ge(p)\leq\JJ\ge[J]$ and $J\ci I$. Since $\vv<\ge(p),\ge[I]>$
is a minimal pair, we conclude that $I = J = \set{x,y}$.
\end{proof}

\begin{claim} For any $p$, $q\in P$,
\[
\ge(p)\DD_L\ge(q)\q \text{if{f}}\q p\dtr q. \]
\end{claim}

\begin{proof} The fact that
\[
\ge(p)\DD_L\ge(q)\q \text{implies that}\q p\dtr q \] follows
immediately from Claim~2.

Conversely, assume that $p\dtr q$. In particular, $X=[p]^{\utr}$
contains
$\set{p,q}$, thus, since $p\ne q$ and $|X|\ne2$, there exists
$x\in P-\set{p,q}$ in
$[p]^{\utr}$. By Claim~2, $\vv<\ge(p),\set{\ge(q),\ge(x)}>$ is a
minimal pair of
$L$; whence $\ge(p)\DD_L\ge(q)$.
\end{proof}

\begin{claim} For any $p$, $q\in P$,
\[
\ge(p)\utr_L\ge(q)\q \text{if{f}}\q p\utr q. \]
\end{claim}

\begin{proof} The fact that
\[
\ge(p)\utr_L\ge(q)\q \text{implies that}\q p\utr q \] follows
immediately from Claim~3 and the fact that the atoms of $L$ are
the $\ge(p)$, for $p\in P$. Conversely, suppose that $p\utr q$.
If $p\ne q$, then, by Claim~3,
$\ge(p)\DD_L\ge(q)$, thus, \emph{a fortiori}, $\ge(p) \utr_L
\ge(q)$. If
$p=q$, then $\ge(p)=\ge(q)$, thus, \emph{a fortiori},
$\ge(p)\utr_L\ge(q)$.
\end{proof}

This last claim concludes the proof of Theorem~\ref{T:Main}.
\end{proof}

\section{Partitions of a finite set}

If $L$ is a finite lattice, then the kernel of the canonical map
$\gQ\colon\J(L)\twoheadrightarrow\J(\Con L)$, namely,
$\asymp_L$ (see Theorem~\ref{T:DCon}), defines a partition of
$\J(L)$. The following result describes exactly what kind of
partition this can be.

\begin{theorem}\label{T:Part} Let $P$ be a finite set, let
$\asymp$ be an equivalence relation on $P$. Then the following
conditions are equivalent:
\begin{enumerate}
\item[\tup{(i)}] There exists a finite atomistic lattice $L$
such that
$\vv<P,\asymp>\iso\vv<\J(L),\asymp_L>$.

\item[\tup{(ii)}] There exists a finite lattice $L$ such that
$\vv<P,\asymp>\iso\vv<\J(L),\asymp_L>$.

\item[\tup{(iii)}] There exists $p\in P$ such that
$|[p]^{\asymp}|\ne2$.
\end{enumerate}
\end{theorem}

\begin{proof}\hfill

(i)$\Rightarrow$(ii) is trivial.

(ii)$\Rightarrow$(iii). Let $L$ be a finite lattice. We prove
that
$\vv<\J(L),\asymp_L>$ satisfies the condition of (iii). Suppose,
to the contrary, that all $\asymp_L$-equivalence classes have
exactly two elements. Let $p\in\J(L)$ be $\utr_L$-maximal, in
the sense that $p\utr_L x$ implies
$p\asymp_Lx$, for all
$x\in\J(L)$. Then $[p]^{\utr_L}=\set{p,p'}$ for the other
element $p'$ of
$[p]^{\asymp_L}$, which contradicts Proposition~\ref{P:Constr}.

(iii)$\Rightarrow$(i). Let $\vv<P,\asymp>$ satisfy (iii). Let
$a\in P$ such that $|[a]^{\asymp}|\ne2$. For $p$, $q\in P$, we
say that $p\utr q$ holds, if either $p\asymp q$, or
$|[p]^{\asymp}|=2$ and $q\asymp a$. It is straightforward to
verify the following statements:
\begin{itemize}
\item[(a)] $\utr$ is a quasi-ordering on $P$. \item[(b)]
$p\asymp q$ if{f}
$p\utr q\utr p$, for any $p$, $q\in P$. \item[(c)]
$|[p]^{\utr}|\ne2$, for any
$p\in P$. \end{itemize}

Therefore, by Theorem~\ref{T:Main}, there exists a finite
atomistic lattice
$L$ such that $\vv<\J(L),\utr_L>\iso\vv<P,\utr>$. By (b) above,
$\vv<\J(L),\asymp_L>\iso\vv<P,\asymp>$.
\end{proof}

\section{The minimal number of
join-irreducibles}\label{S:MinNo}

Let $D$ be a finite distributive lattice. In this section, we
shall compute the minimal number of join-irreducible elements in
a finite lattice $L$ such that $\Con L$ is isomorphic to $D$. If
$P$ is the poset of join-irreducible elements of $D$, then $D$
is isomorphic to $\HH(P)$, the poset of hereditary subsets of
$P$, which makes it possible to formulate the problem in terms
of the finite poset $P$. We shall first assign to $P$ a natural
number $\ga(P)$.

\begin{definition}\label{D:spike} Let $P$ be a poset. A
\emph{spike} of $P$ is a pair $\vv<p,q>$ of elements of $P$ such
that $q$ is maximal and $q$ is the unique element of $P$ that
covers $p$. We define
\begin{align*}
\pp P&=\setm{q\in P}{\vv<p,q>\text{ is a spike of }P, \text{ for
some }p\in P};\\
\ppu P&=\setm{q\in P}{\vv<p,q>\text{ is a spike of }P, \text{
for a unique }p\in P};\\
\ppm P&=\pp P-\ppu P;\\
\ga(P)&=|\ppu P|+2|\ppm P|.
\end{align*}
\end{definition}

In particular, we say that $P$ is \emph{spike-free}, if there
are no spikes in
$P$. Note that $P$ is spike-free if{f} $\ga(P)=0$. Equivalently,
$|\,[p]^{\leq}\,|\ne2$, for any $p\in P$; note at this point the
similarity with the condition in Proposition~\ref{P:Constr}.

\begin{theorem}\label{T:Ineq} Let $L$ be a finite lattice. Then
the following inequality holds: \[ |\J(L)|\geq|\J(\Con
L)|+\ga(\J(\Con L)). \]
\end{theorem}

\begin{proof} Let $P=\J(L)$ and $\ol{P}=\J(\Con L)$. For any
$p\in P$, as in the statement of Theorem~\ref{T:DCon}, denote by
$\gQ(p)$ the principal congruence of
$L$ generated by the pair $\vv<p_*,p>$.

\setcounter{claim}{0}
\begin{claim} Let $\vv<\ol{p},\ol{q}>$ be a spike of $\ol{P}$.
Then either
$|\gQ^{-1}\set{\ol{p}}|\geq2$ or $|\gQ^{-1}\set{\ol{q}}|\geq3$.
\end{claim}

\begin{proof} Assume that the conclusion of the claim does not
hold. Then
$\gQ^{-1}\set{\ol{p}}$ is a singleton, say, $\set{p}$, and
$\gQ^{-1}\set{\ol{q}} = \set{q,q'}$, for some
$q$, $q'\in P$. If $q=q'$, then $[p]^{\utr_L}=\set{p,q}$ has
exactly two elements, a contradiction by
Proposition~\ref{P:Constr}. So $q\ne q'$. Since
$\gQ(q) = \gQ(q') = \ol{q}$, the set $[q]^{\utr_L}=\set{q,q'}$
has, again, exactly two elements, a contradiction by
Proposition~\ref{P:Constr}.
\end{proof}

For all $\ol{q}\in\pp\ol{P}$, we define $X(\ol{q}) \ci \ol{P}$,
by
\[ X(\ol{q})=\set{\ol{q}}\uu
\setm{\ol{p}\in\ol{P}}{\vv<\ol{p},\ol{q}>\text{ is a spike of
}\ol{P}}. \] By the definition of a spike, the sets $X(\ol{q})$,
for $\ol{q}\in\ol{P}$, are mutually disjoint. Furthermore, it
follows immediately from Claim~1 that for all $\ol{q}\in\ol{P}$,
the following statements hold: \begin{align*}
|\,\gQ^{-1}[X(\ol{q})]\,|\geq|X(\ol{q})|+1,\q &\text{for }
\ol{q}\in\ppu\ol{P};\\
|\,\gQ^{-1}[X(\ol{q})]\,|\geq|X(\ol{q})|+2,\q &\text{for }
\ol{q}\in\ppm\ol{P}.
\end{align*} Let $X=\UUm{X(\ol{q})}{\ol{q}\in\ol{P}}$; then \[
|\,\gQ^{-1}[X]\,| \geq|X| + \ga(\ol{P}). \] Since the map $\gQ$
is surjective,
\[ |P|\geq|\ol{P}|+\ga(\ol{P}),
\] which is the desired conclusion.
\end{proof}

The converse of Theorem~\ref{T:Ineq} is provided by the
following result, which proves that the bound $\ga(\J(\Con L))$
is best possible:

\begin{theorem}\label{T:OptJoin} Let $P$ be a finite poset. Then
there exists a finite atomistic lattice $L$ such that $\Con
L\iso\HH(P)$ and
$|\J(L)|=|P|+\ga(P)$. \end{theorem}

\begin{proof} Define $P_1 \ci P$ as follows:
\[ P_1=\setm{p\in P}{\vv<p,q>\text{ is a spike, for some
}q\in\ppu P}. \] Note that for $p\in P_1$, there exists a
\emph{unique} $q\in\ppu P$ that covers
$p$. In particular, $|P_1|=|\ppu P|$.

Then we define a finite set $Q$, by
\[ Q=(P-(P_1\uu\ppm P))\uu(P_1\times 2)\uu(\ppm P\times 3) \] (a
disjoint union), where we identify $2$ with $\set{0,1}$ and $3$
with $\set{0,1,2}$. Note that
\[ |Q|=|P|+|P_1|+2|\ppm P|=|P|+\ga(P).
\] Let $\gp\colon Q\twoheadrightarrow P$ be the natural
projection, that is,
$\gp(x)=x$, if $x\in P-(P_1\uu\ppm P)$, and $\gp(\vv<x,i>)=x$,
if $\vv<x,i>\in P\times 3$. We define a quasi-ordering $\utr$ on
$Q$, by
\[ p\utr q\q\text{if{}f}\q \gp(p)\leq\gp(q). \] We now verify
that $\utr$ satisfies the assumption (iii) of
Theorem~\ref{T:Main}. So let $p\in Q$; we shall prove that
$[p]^{\utr}$ does not have exactly two elements. Let us assume
otherwise, that is, let
\begin{equation}\label{Eq:Fin2}
[p]^{\utr}=\set{p,q},\quad\text{for some }q\in Q-\set{p}.
\end{equation} We separate three cases.

\setcounter{case}{0}
\begin{case}
$p=\vv<x,i>$, where $x\in P_1$ and $i<2$. \end{case} By the
definition of
$P_1$, there exists $y\in\ppu P$ such that $\vv<x,y>$ is a spike
of $P$. Note that $y$ belongs to $P-(P_1\uu\ppm P)$, thus to
$Q$, so that $y$ belongs to
$[p]^{\utr}$. Therefore,
\[ [p]^{\utr}\text{ contains }\set{\vv<x,0>,\vv<x,1>,y}, \]
which contradicts
\eqref{Eq:Fin2}.

\begin{case}
$p=\vv<x,i>$, where $x\in\ppm P$ and $i<3$. \end{case} Then
$[p]^{\utr}$ equals $\set{\vv<x,0>,\vv<x,1>,\vv<x,2>}$, which
contradicts \eqref{Eq:Fin2} again.

\begin{case}
$p\in P-(P_1\uu\ppm P)$.
\end{case} If $p$ is maximal in $P$, then $\gp(q)=\gp(p)=p$
belongs to
$P-(P_1\uu\ppm P)$, thus $q$ belongs to $P-(P_1\uu\ppm P)$, so
that $p=q$, which contradicts \eqref{Eq:Fin2}. Hence $p$ is not
maximal in $P$. If $p$ is not the bottom element of a spike in
$P$, then there are distinct $x$, $y$ in
$P$ such that $p<x$, $y$. If
$x'$, $y'\in Q$ are such that $\gp(x')=x$ and $\gp(y')=y$, then
$[p]^{\utr}$ contains the three-element set $\set{p,x',y'}$,
which contradicts
\eqref{Eq:Fin2}. So there exists $r\in P$ such that $\vv<p,r>$
is a spike of
$P$. Since $p$ does not belong to $P_1$, $r$ belongs to $\ppm
P$. Hence
$[p]^{\utr}$ contains the four-element set
$\set{p,\vv<r,0>,\vv<r,1>,\vv<r,2>}$, which contradicts
\eqref{Eq:Fin2}.

By Theorem~\ref{T:Main}, there exists a finite atomistic lattice
$L$ such that
$\vv<\J(L),\utr_L>\iso\vv<Q,\utr>$. In particular,
\[ |\J(L)|=|Q|=|P|+\ga(P).
\] Furthermore, $\J(\Con L)$ is isomorphic to the quotient of
$\vv<\J(L),\utr_L>$ by the equivalence relation associated with
$\utr_L$, thus to the quotient of $\vv<Q,\utr>$ by the
equivalence relation associated with
$\utr$. Since the latter is exactly the kernel of $\gp$, the
corresponding quotient is isomorphic to $\vv<P,\leq>$. Hence
\[
\J(\Con L)\iso P
\] (as posets), from which it follows that $\Con L\iso\HH(P)$.
\end{proof}

As an immediate consequence of Theorem~\ref{T:OptJoin}, we
obtain the following result:

\begin{corollary}\label{C:5/3bound} Let $D$ be a finite
distributive lattice. Then there exists a finite atomistic
lattice $L$ such that $\Con L\iso D$ and
\begin{equation}\label{Eq:5/3bound}
|\J(D)|\leq|\J(L)|\leq\frac{5}{3}|\J(D)|.
\end{equation} Furthermore, the constant $5/3$ in the inequality
\eqref{Eq:5/3bound} is best possible.
\end{corollary}

\begin{proof} Put $P=\J(D)$. By Theorem~\ref{T:OptJoin}, to
establish the inequality \eqref{Eq:5/3bound}, it suffices to
establish the inequality
\begin{equation}\label{Eq:Alph2/3}
\ga(P)\leq\frac{2}{3}|P|.
\end{equation} We put
\begin{align*} X(q)&=\set{q}\uu\setm{p\in P}{\vv<p,q>\text{ is
a spike of }P},
\q\text{for all }q\in\pp P,\\ P_\mathrm{u}&=\UUm{X(q)}{q\in\ppu
P},\\ P_\mathrm{m}&=\UUm{X(q)}{q\in\ppm P}.
\end{align*}
As in the proof of
Theorem~\ref{T:Ineq}, we note that the sets $X(q)$, for $q\in\pp
P$, are mutually disjoint. Furthermore, for
$q\in\ppu P$, $|X(q)|=2$, while for $q\in\ppm P$, $|X(q)|\geq
3$. It follows that $|P_\mathrm{u}|=2|\ppu P|$ and
$|P_\mathrm{m}|\geq 3|\ppm P|$. Therefore,
\begin{align*}
\ga(P)&=|\ppu P|+2|\ppm P|\\
&\leq\frac{1}{2}|P_\mathrm{u}|+\frac{2}{3}|P_\mathrm{m}|\\
&\leq\frac{2}{3}|P|,
\end{align*} which completes the proof of \eqref{Eq:Alph2/3}.

The upper bound $\frac{5}{3}|\J(D)|$ in the inequality
\eqref{Eq:5/3bound} is reached by defining $P$ as the
three-element set $\set{u,v,1}$, endowed with the ordering
defined by $u$, $v<1$. For this example,
$\ga(P)=2=\frac{2}{3}|P|$.
\end{proof}

\begin{corollary}\label{C:SpFree} Let $P$ be a finite poset.
Then the following are equivalent:
\begin{enumerate}
\item[\tup{(i)}] There exists a finite, atomistic, lower bounded
lattice $L$ such that
\[
\Con L\iso\HH(P).
\]
\item[\tup{(ii)}] There exists a finite lower bounded lattice
$L$ such that
\[
\Con L\iso\HH(P).
\]
\item[\tup{(iii)}] $P$ is spike-free.
\end{enumerate}
\end{corollary}

\smallskip

\tbf{Note}
Lower bounded finite lattices were introduced in
R.\,N. McKenzie
\cite{rM70}, see also R.~Freese, J.~{}Je\v zek, and J.\,B. Nation
\cite{FJN95}. A finite lattice $A$ is lower bounded if{}f $A$
has no
$\DD_A$-cycle. An equivalent condition is that $A$ be
\emph{amenable}: The tensor product $A \otimes L$ is a lattice,
for every lattice $L$ with $0$; see G. Gr\"atzer and F. Wehrung
\cite{GW0,GW1}.
\smallskip

\begin{proof}\hfill

(i)$\Rightarrow$(ii) is trivial.

If $L$ is lower bounded, then $|\J(L)|=|\J(\Con L)|$, see
Lemma~2.40 in
\cite{FJN95}, thus, by Theorem~\ref{T:Ineq}, $\ga(\J(\Con
L))=0$, that is,
$\J(\Con L)$ is spike-free. This proves that
(ii)$\Rightarrow$(iii).

If $P$ is spike-free, that is, $\ga(P)=0$, then, by
Theorem~\ref{T:OptJoin}, there exists a finite atomistic lattice
$L$ such that $|\J(L)|=|P|$ and $\Con L\iso\HH(P)$. {}From the
second equality it follows that $\J(\Con L)\iso P$, whence
$|\J(L)\,|=|\J(\Con L)\,|$. Again by Lemma~2.40 in
\cite{FJN95}, $L$ is lower bounded.
\end{proof}

The following example shows that in (i) of
Corollary~\ref{C:SpFree}, one cannot replace ``atomistic'' by
the stronger condition ``sectionally complemented''.

\begin{example}\label{Ex:nosc} A finite, spike-free poset $P$
such that there exists no finite, lower bounded,
\emph{sectionally complemented} lattice $L$ such that $\Con
L\iso\HH(P)$.
\end{example}

\begin{proof} Let $P=\set{p,q,q_0,q_1}$, and let the ordering of
$P$ be generated by the pairs
$p<q$, $q<q_0$, and $q<q_1$. It is obvious that $P$ is
spike-free. Assume that
$P\iso\Con L$ for some finite, lower bounded, sectionally
complemented lattice
$L$. Note, in particular, that $L$ is atomistic. Since $L$ is
lower bounded and finite, $\vv<\J(L),\utr_L>$ is isomorphic to
$\vv<P,\leq_P>$. Thus, without loss of generality,
$\vv<\J(L),\utr_L>=\vv<P,\leq_P>$.

In particular, $q\tr_Lq_0$, thus there exists $x\in P$ such that
$q\DD_Lx$ and
$x\utr_Lq_0$. The first condition implies that
$x\in\set{q_0,q_1}$, and the second condition implies then that
$x=q_0$; whence $q\DD_Lq_0$. By Lemma~\ref{L:DescrD}, there
exists a subset $I$ of $P$ such that $\vv<q,I>$ is a minimal
pair of $L$ and
$q_0\in I$. For all $x\in I$, $q\utr_Lx$ and $q\ne x$, so we
obtain that
$x\in\set{q_0,q_1}$. Thus $I\ci\set{q_0,q_1}$. Since
$|I|\geq2$, it follows that
$I=\set{q_0,q_1}$. In particular, we obtain the inequality
\begin{equation}\label{Eq:qq0q1} q<q_0\jj q_1.
\end{equation}

Since $p\tr_Lq$, there exists $J\ci P$ such that $\vv<p,J>$ is a
minimal pair of
$L$. {}From $J\ci P-\set{p}$, it follows that $p<q\jj q_0\jj
q_1$. Thus, by
\eqref{Eq:qq0q1}, $p<q_0\jj q_1$ and
\begin{equation}\label{Eq:q0q1x} q_0\jj q_1=1.
\end{equation}

Now let $x$ be a complement of $q$ in $L$. Without loss of
generality, $q_0\ne x$. Furthermore, note that $q_0\leq 1=x\jj
q$. Thus $x$ cannot be an atom of
$L$; otherwise, since $q_0\ne x$, we have $q_0\DD_Lx$, which is
impossible. Moreover, $x\ne q_0\jj q_1$ by \eqref{Eq:q0q1x}.
Since $x$ is a join of atoms distinct from $q$, it follows that
$x=p\jj q_i$, for some $i<2$. Therefore,
\[ q_{1-i}\leq 1=q\jj x=p\jj q\jj q_i.
\] Since $p$, $q$, $q_0$, and $q_1$ are atoms of $P$, there
exists
$X\ci\set{p,q,q_i}$ such that $\vv<p,X>$ is a minimal pair of
$L$. So
$q_{1-i}\DD_Ly$, for all $y\in X$, thus for some
$y\in\set{p,q,q_i}$, a contradiction.
\end{proof}

\section{Open problems}

\begin{problem} Characterize the relation $\utr_L$, for a finite
sectionally complemented lattice
$L$.
\end{problem}

By Example~\ref{Ex:nosc}, not every relation of the form
$\utr_K$, for $K$ finite and atomistic, is of the form $\utr_L$,
for $L$ finite and sectionally complemented.

\begin{problem} Let $D$ be a finite distributive lattice. Find a
simple way of computing the least possible value of $|\J(L)|$,
for a finite, sectionally complemented lattice $L$ such that
$\Con L\iso D$.
\end{problem}

By Example~\ref{Ex:nosc}, the least possible value of $|\J(L)|$,
for a finite, sectionally complemented lattice $L$ such that
$\Con L\iso D$, may be larger than the least possible value of
$|\J(L)|$, for a finite, atomistic lattice
$L$ such that $\Con L\iso D$.

A related question is the following:

\begin{problem}\label{Pb:BestCstSC} Determine the least
constant $k$ such that for every finite distributive lattice
$D$, there exists a finite, sectionally complemented lattice $L$
such that $\Con L\iso D$ and
\[ |\J(L)|\leq k|\J(D)|.
\]
\end{problem}

By \cite{GS62}, $k$ is less or equal than $2$. The value of the
constant defined similarly for the class of atomistic lattices
(or the class of all lattices as well) equals $5/3$, by
Corollary~\ref{C:5/3bound}.

\begin{problem} Characterize the relation $\DD_L$ for a finite
(resp., finite atomistic, finite sectionally complemented)
lattice $L$.
\end{problem}

\begin{problem} Let $\E V$ be a variety of lattices. If $D$ is a
finite distributive lattice representable by a finite lattice in
$\E V$, compute the least possible value of $|\J(L)|$, for a
finite lattice $L$ in $\E V$ such that $\Con L\iso D$.
\end{problem}

\begin{problem} Let $\asymp$ be an equivalence relation on a
finite set $P$ such that $|[a]^{\asymp}|\ne2$, for some $a\in
P$. Does there exists a finite, sectionally complemented lattice
$L$ such that
$\vv<\J(L),\asymp_L>\iso\vv<P,\asymp>$?
\end{problem}

In the proof of Theorem~\ref{T:Part}, we construct a finite
atomistic lattice
$L$ such that $\vv<\J(L),\asymp_L>\iso\vv<P,\asymp>$, however,
this lattice
$L$ may not be sectionally complemented.


\begin{thebibliography}{9}

\bibitem{hG73} H.\,S. Gaskill,
\emph{On transferable semilattices}, Algebra
Universalis~\tbf{2} (1973), 303--316.

\bibitem{GGP75} H.\,S. Gaskill, G. Gr\"atzer, and C.\,R. Platt,
\emph{Sharply transferable lattices}, Canad. J. Math.~\tbf{28}
(1975), 1246--1262.

\bibitem{FJN95} R.~Freese, J.~{}Je\v zek, and J.\,B. Nation,
\emph{Free lattices}, Mathematical Surveys and Monographs,
Vol.~\tbf{42}, American Mathematical Society, Providence, RI,
1995. viii+293 pp.

\bibitem{GS62} G. Gr\"atzer and E.\,T. Schmidt,
\emph{On congruence lattices of lattices}, Acta Math. Acad. Sci.
Hungar.~\tbf{13} (1962), 179--185.

\bibitem{GS}
\bysame,
\emph{Congruence-preserving extensions of finite lattices into
sectionally complemented lattices}, Proc. Amer. Math. Soc.,
\tbf{127} (1999), 1903--1915.

\bibitem{GW0} G. Gr\"atzer and F. Wehrung,
\emph{Tensor products of lattices with zero, revisited}, J. Pure
Appl. Algebra
\tbf{147} (2000), 273--301.

\bibitem{GW1}
\bysame,
\emph{Tensor products and transferability of semilattices,}
Canad. J. Math.
\tbf{51} (1999), 792--815.

\bibitem{rM70} R.\,N. McKenzie,
\emph{Equational bases and nonmodular lattice varieties,}
Trans.\ Amer.\ Math.\ Soc.\ \tbf{174} (1972), 1--43.

\bibitem{Tisch} M. Tischendorf,
\emph{The representation problem for algebraic distributive
lattices}, Ph. D. Thesis, TH Darmstadt, 1992.

\end{thebibliography}
\end{document}